\documentclass[11pt, twoside]{article}
\usepackage{ amssymb, amsmath,amsthm}
\usepackage{latexsym}
\usepackage{float,latexsym,shortvrb}
\usepackage[dvips]{graphics, epsfig, color}

\newcommand{\ignore}[1]{}
\ignore{The text written between the parentheses here is ignored by latex}

\def\@begintheorem#1#2{\par\bgroup{\sc #1\ #2. }\it\ignorespaces}
\def\@opargbegintheorem#1#2#3{\par\bgroup{\sc #1\ #2\ (#3). }\it\ignorespaces}
\def\@endtheorem{\egroup}

\newtheorem{theorem}{Theorem}[section]
\newtheorem{corollary}[theorem]{Corollary}
\newtheorem{lemma}[theorem]{Lemma}
\newtheorem{remark}[theorem]{Remark}
\newtheorem{proposition}[theorem]{Proposition}
\newtheorem{definition}[theorem]{Definition}
\newtheorem{example}[theorem]{Example}
\newtheorem{question}[theorem]{Question}
\newtheorem{problem}[theorem]{Problem}
\newtheorem{conjecture}[theorem]{Conjecture}

\newcommand{\bd}[1]{\begin{definition}\rm\label{#1}}
\newcommand{\bt}[1]{\begin{theorem}\label{#1}}
\newcommand{\bc}[1]{\begin{corollary}\label{#1}}
\newcommand{\bcj}[1]{\begin{conjecture}\label{#1}}
\newcommand{\bl}[1]{\begin{lemma}\label{#1}}
\newcommand{\bp}[1]{\begin{proposition}\label{#1}}
\newcommand{\be}[1]{\begin{example}\rm\label{#1}}
\newcommand{\bq}[1]{\begin{question}\rm\label{#1}}
\newcommand{\bprob}[1]{\begin{problem}\rm\label{#1}}
\newcommand{\beq}[1]{\begin{eqnarray}\label{#1}}
\newcommand{\br}[1]{\begin{remark}\rm\label{#1}}
\newcommand{\bpr}{\par{\it Proof}. \ignorespaces}

\newcommand{\el}{\end{lemma}}
\newcommand{\ep}{\end{proposition}}
\newcommand{\ee}{\end{example}}
\newcommand{\eq}{\end{question}}
\newcommand{\eprob}{\end{problem}}
\newcommand{\ecj}{\end{conjecture}}
\newcommand{\eeq}{\end{eqnarray}}
\newcommand{\ed}{\end{definition}}
\newcommand{\et}{\end{theorem}}
\newcommand{\ec}{\end{corollary}}
\newcommand{\er}{\end{remark}}
\newcommand{\epr}{{\ \vbox{\hrule\hbox{%
   \vrule height1.3ex\hskip0.8ex\vrule}\hrule}}\\\par}

\begin{document}

 \pagestyle{myheadings}
 \markboth{Berge sorting}{Antoine Deza and William Hua}

\thispagestyle{empty}
\title{{\bf Berge Sorting}\\
 }
\author{Antoine  Deza
\and
William Hua}
\date{July 14, 2005}
\maketitle
\begin{abstract}
In 1966, Claude Berge proposed the following sorting problem. 
Given a string of $n$ alternating
white and black pegs on a one-dimensional board consisting of an unlimited number of empty
holes, rearrange the pegs into a string consisting of $\lceil\frac{n}{2}\rceil$
white pegs followed immediately  by $\lfloor\frac{n}{2}\rfloor$ black pegs 
(or vice versa) using only moves which take $2$ adjacent pegs to $2$ vacant
adjacent holes. Avis and Deza proved that the alternating string can be 
sorted in $\lceil\frac{n}{2}\rceil$ such {\em Berge $2$-moves} for $n\geq 5$. 
Extending Berge's original problem, we consider the same sorting problem using  {\em Berge $k$-moves}, i.e., moves
 which take $k$ adjacent pegs to $k$ vacant adjacent holes. We prove that the alternating string 
can be sorted in $\lceil\frac{n}{2}\rceil$  Berge $3$-moves for $n\not\equiv 0\pmod{4}$ 
and in $\lceil\frac{n}{2}\rceil+1$  Berge $3$-moves for $n\equiv 0\pmod{4}$, for $n\geq 5$.
In general, we conjecture that, for any $k$ and large enough $n$, the alternating string can be sorted in 
$\lceil\frac{n}{2}\rceil$  Berge $k$-moves.
This estimate is tight as  $\lceil\frac{n}{2}\rceil$ is a lower bound for the minimum number of required Berge $k$-moves for $k\geq 2$
and $n\geq 5$.
\end{abstract}

\section{Introduction}
In a column 
that appeared in the Revue Fran\c{c}aise de Recherche Op\'erationnelle 
in 1966, 
entitled {\em  Probl\a`emes plaisans et d\a'electables}
in homage to the 17th century work of Bachet~\cite{Ba12},
Claude Berge~\cite{Be66} proposed the following sorting problem:

\begin{quote}
For $n\geq 5$, given a string of $n$ alternating
white and black pegs on a one-dimensional board consisting of an unlimited number of empty
holes, we are required to rearrange the pegs into a string consisting of $\lceil\frac{n}{2}\rceil$
white pegs followed immediately  by $\lfloor\frac{n}{2}\rfloor$
 black pegs (or vice versa) using only moves which take $2$ adjacent pegs to $2$ vacant
adjacent holes. Berge noted that the minimum number of moves required is $3$ for $n=5$ and $6$,
 and $4$ for $n=7$. See Figure~\ref{5-2} for a sorting of $5$ pegs in $3$ moves.
\begin{figure}[hbt]
\begin{center}
\epsfig{file=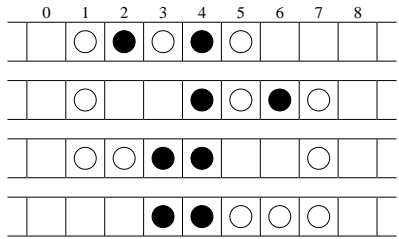,width=5.2cm}
\caption{Sorting $5$ pegs in $3$ moves}
\label{5-2}
\end{center}
\end{figure}
\end{quote}

\newpage
Avis and Deza~\cite{AD0X} provided a solution in  $\lceil\frac{n}{2}\rceil$ {\em Berge $2$-moves} for $n\geq 5$. 
Extending Berge's original problem, we consider the same sorting question using only Berge $k$-moves, i.e., moves
 which take $k$ adjacent pegs to $k$ vacant adjacent holes. 
We provide a solution  
in $\lceil\frac{n}{2}\rceil$ Berge $3$-moves for $n\not\equiv 0\pmod{4}$ 
and in $\lceil\frac{n}{2}\rceil+1$  Berge $3$-moves for $n\equiv 0\pmod{4}$  and $n\geq 5$.
The authors generated minimal solutions by computer for a large number of $k$ and $n$ which turned out all be equal to $\lceil\frac{n}{2}\rceil$
except for the few first small values of $n$. 
Note that, for $k\geq 2$, $\lceil\frac{n}{2}\rceil$ is a lower bound for the 
minimum number of required Berge $k$-moves, see Section~\ref{mini}.
To the best of our knowledge, this property was not noticed earlier. We
conjecture that for any  $k$ and large enough $n$,
the alternating string 
can be rearranged into a string consisting of $\lceil\frac{n}{2}\rceil$
white pegs followed immediately  by $\lfloor\frac{n}{2}\rfloor$
black pegs (or vice versa) by only $\lceil\frac{n}{2}\rceil$
moves which take $k$ adjacent pegs to $k$ vacant adjacent holes.

\section{Notation}
We follow and adapt the notation used in~\cite{AD0X,Be66}. The starting game board consists
of  $n$ alternating white and black pegs sitting in the positions $1$ through $n$. 
A single Berge $k$-move will be denoted as $\{\ j\ i\ \}$, in which case, the pegs in the positions $i,i+1,\dots,i+k-1$ 
are moved to the vacant holes $j, j+1,\dots,j+ k-1$. Successive moves are concatenated as $\{\ j\ i\ \}\ \cup\ \{\ l\ k\ \}$, 
which means perform $\{\ j\ i\ \}$ followed by $\{\ l\ k\ \}$.  Often, a move  fills an empty hole created as an effect of the previous move, 
and the 
resulting notation $\{\ j\ k\ \}\ \cup\ \{\ k\ i\ \}$ is abbreviated as $\{\ j\ k\ i\ \}$. This can be extended to more than two such moves as well.
${\cal S}_{n,k}$ denotes a solution for $n$ pegs by Berge $k$-moves and $h(n,k)$ denotes the minimum number 
of required $k$-moves, i.e., the length of a shortest solution. For example, with this notation, the values
$h(5,2)=h(6,2)=3$ and $h(7,2)=4$ given by Berge~\cite{Be66}  are illustrated in Table~\ref{Bergen-2}.
\begin{table}[htb]
\caption{First solutions using Berge $2$-moves}\label{Bergen-2}
\begin{center}
\begin{tabular}{r c l l l l l l l l l l l l l l l l l l l}
${\cal S}_{5,2}$  & =  & \{   6 \  2 \  5 \  1  \}  \\
${\cal S}_{6,2}$ & =  & \{    7 \  4 \ 1  \}  $\cup$  \{   9 \  3   \}   \\
${\cal S}_{7,2}$ & =  & \{    8 \  2 \  5 \ 8 \  1  \}
\end{tabular}
\end{center}
\end{table}

\section{Main Results}
\subsection{Minimum number of required Berge $k$-moves} \label{mini} 
Let  ${\cal D}_{n,k}(i)$ denote
the {\em disorder}, i.e.,  the number of pegs whose right neighbour is not a peg of the same colour after the $i$-th 
Berge $k$-move. One can easily check that $|{\cal D}_{n,k}(i)-{\cal D}_{n,k}(i+1)|\leq 2$. A move such that 
${\cal D}_{n,k}(i)-{\cal D}_{n,k}(i+1)=2$ (resp. $1$ and $0$) is called {\em optimal} (resp. {\em suboptimal} and {\em neutral}).

\bl{k>0}
At least $\lfloor\frac{n}{2}\rfloor$ Berge $k$-moves are required to sort a string of  $n$ alternating white and black pegs.
In other words, $h(n,k)\geq\lfloor\frac{n}{2}\rfloor$ for $k\geq 1$ and $n\geq 3$.
\bpr
The disorder of the initial board is ${\cal D}_{n,k}(0)=n$ and the disorder of the sorted string is ${\cal D}_{n,k}(h(n,k))=2$.
Since the first move cannot be optimal, i.e., ${\cal D}_{n,k}(0)-{\cal D}_{n,k}(1)\leq 1$, and the following moves satisfy
${\cal D}_{n,k}(i)-{\cal D}_{n,k}(i+1)\leq 2$, we have $h(n,k)\geq\lfloor\frac{n}{2}\rfloor$.
\epr 
\el
\begin{table}[htb]
\caption{Sorting $n$ pegs in  $\lfloor\frac{n}{2}\rfloor$ Berge $1$-moves for $n\equiv3\pmod{4}$}\label{Bergen-1}
\begin{center}
\begin{tabular}{r c l l l l l l l l l l l l l l l l l l l}
${\cal S}_{3,1}$ & =  & \{  4  \  1   \}   \\
${\cal S}_{7,1}$ & =  & \{  8 \  3 \ 6  \  1   \}   \\
${\cal S}_{11,1}$ & = & \{  12 \  3 \  10 \ 5 \  8 \ 1   \}  \\ 
${\cal S}_{15,1}$ & = & \{  16 \  3 \  14 \ 5 \  12 \ 7 \ 10  \ 1   \} \\ 
${\cal S}_{4i+3,1}$ & = & \{  4$i$+4 \  3 \  4$i$+2 \ 5 \  4$i$ \ 7 \ 4$i$-2 \ 9 \ \dots \ 2$i$+4  \ 1  \} 
\end{tabular}
\end{center}
\end{table}
Lemma~\ref{k>0}  is tight because, for $k=1$, we have $h(n,1)=\lfloor\frac{n}{2}\rfloor$ for $n\equiv3\pmod{4}$, see Table~\ref{Bergen-1}.
Solutions in $\lceil\frac{n}{2}\rceil$ Berge $1$-moves for  $n\not\equiv3\pmod4$  are very similar to the ones  in  $\lfloor\frac{n}{2}\rfloor$ $1$-moves for
$n\equiv3\pmod{4}$.
Avis and Deza noticed in~\cite{AD0X} that $h(n,2)\geq \lceil\frac{n}{2}\rceil$ for $n\geq 5$. For $k\geq 2$, Lemma ~\ref{k>0} can be strengthen to the following lemma.
\bl{k>1}
At least $\lceil\frac{n}{2}\rceil$ Berge $k$-moves are required to sort a string of  $n$ alternating white and black pegs.
In other words, $h(n,k)\geq\lceil\frac{n}{2}\rceil$ for $k\geq 2$ and $n\geq 5$.
\bpr
As Lemma~\ref{k>0} and~\ref{k>1} are equivalent for even $n$, let us assume that, for odd $n\geq 5$, we
have a solution in $\lfloor\frac{n}{2}\rfloor$ Berge $k$-moves. 
It implies that, after the first suboptimal move, all the following moves are optimal. 
We derive a contradiction for $k=3$ and the same argument can be used for any $k\geq 2$.
Since $n$ is odd, the initial board is something like 
$\circ\bullet\circ\bullet\circ\bullet\circ\bullet\circ\bullet\circ$ where $\circ$ and $\bullet$ represent
white and black pegs. By symmetry, we can assume the first move is to the right.  This first suboptimal move has 
to take $3$ pegs from the interior of the string to the position $n+1$. For example, with $n=11$, the board
after the first move is something like $\circ\bullet - - - \bullet\circ\bullet\circ\bullet\circ\circ\bullet\circ$.  
The next move must fill the vacancy
with a $\bullet\star\bullet$ triple, where $\star$ is any colour, but additionally the $\bullet\star\bullet$ triple must have 
been taken from between two white pegs to maintain optimality.  Similarly,
the subsequent moves must alternate between optimal fillings of $\bullet - -  -\bullet$ and $\circ - -  - \circ$ vacancies. 
Consider the last $4$ (or $k+1$ in general) pegs, $\circ\circ\bullet\circ$, after the first suboptimal move: 
As the last triple, $\circ\bullet\circ$, or the triple before, $\circ\circ\bullet$, do not correspond to an optimal filling,
the black (or white) peg  in the last 2 positions cannot be sorted by optimal moves.
\epr 
\el

\subsection{Optimal solutions for sorting by Berge $k$-moves}
We first recall that a solution for sorting the alternating string in $\lceil\frac{n}{2}\rceil$ Berge $2$-moves for $n\geq 5$
was given in~\cite{AD0X}.
\bp{2moves}
{\em \cite{AD0X}}
For $n\geq 5$, a string of  $n$ alternating white and black pegs can be sorted  
in $\lceil\frac{n}{2}\rceil$ Berge $2$-moves.
In other words, $h(n,2)=\lceil\frac{n}{2}\rceil$ for  $n\geq 5$.
\ep
Considering the case $k=3$, we prove that $h(n,3)= \lceil\frac{n}{2}\rceil$ for $n\not\equiv 0\pmod 4$
and, while computer calculations and preliminary attempts strongly indicated that the same holds for $n\equiv 0\pmod 4$
and $n\geq 20$,
so far we could only exhibit a solution in $\lceil\frac{n}{2}\rceil+1$ Berge $3$-moves for $n\equiv 0\pmod 4$ and $n \geq 8$.

\bp{3moves}
For $n\geq 5$, a string of  $n$ alternating white and black pegs can be sorted  
in $\lceil\frac{n}{2}\rceil$ Berge $3$-moves
for  $n\not\equiv 0\pmod{4}$ 
and in $\lceil\frac{n}{2}\rceil+1$  Berge $3$-moves for $n\equiv 0\pmod{4}$.
In other words, for $n\geq 5$, $h(n,3)=\lceil\frac{n}{2}\rceil$ for $n\not\equiv 0\pmod{4}$ and 
$\lceil\frac{n}{2}\rceil\leq h(n,3)\leq \lceil\frac{n}{2}\rceil+1$ for $n\equiv 0\pmod{4}$. 
\bpr
See Section~\ref{proofs-3} for a description of the solutions ${\cal S}_{n,3}$.
\epr
\ep
Propositions~\ref{2moves} and~\ref{3moves}  
lead to the following conjecture.
\bcj{k_Berge}
For any $k$, a string of  $n$ alternating white and black pegs can be sorted  
in $\lceil\frac{n}{2}\rceil$ Berge $k$-moves for $n\geq  2k+11$.
In other words, $h(n,k)=\lceil\frac{n}{2}\rceil$ for  $k\geq 2$ and $n\geq 2k+11$.
\ecj
To substantiate  Conjecture~\ref{k_Berge}, the authors calculated the values of $h(n,k)$ by computer for  
$k\leq 14$ and $n\leq 50$  and, for these preliminary computations, did not find any counterexample.
See Table~\ref{substance}, which gives the values of $h(n,k)-\lceil\frac{n}{2}\rceil$  for $k \leq 14$ and $n \leq 50$.
Note that the alternating string obviously cannot be sorted by any number of $k$-moves for $n\leq k+1$.
The more conservative conjecture consisting in replacing ^^ ^^ {\em $n\geq 2k+11$}" by  ^^ ^^ {\em  $n\geq {k+2\choose 2}+7$}"
is also consistent with the computations reported in Table~\ref{substance}.
See~\cite{DH05} for detailed and updated computational results.

\bp{any_k}
Let $n\geq k+2\geq 6$, if the following conditions $(i)$, $(ii)$, $(iii)$ and
$(iv)$ are satisfied
\begin{itemize}
\item[$(i)$]
$h(4k,k)=2k$,
\item[$(ii)$] the solution ${\cal S}_{4k,k}$ shifts the string $k$ spaces to the right overall and
is made of moves involving only positions greater than or equal to $1$, 
\item[$(iii)$] for a given $n_k$,
$h(n,k)=\lceil\frac{n}{2}\rceil$ for $n_k\leq n<n_k+4k$,
\item[$(iv)$]
 the $4k$ consecutive solutions ${\cal S}_{n,k}$ for 
$n_k\leq n<n_k+4k$ 
are made of moves involving only positions greater than or equal to $1-k$, 
\end{itemize}
then we have $h(n,k)\leq \lceil\frac{n}{2}\rceil+ \lceil\frac{n}{k}\rceil-\lceil\frac{n_k}{k}\rceil$ for $n\geq n_k$. 
For example, one can check that the conditions $(i)$, $(ii)$, $(iii)$ and $(iv)$ are satisfied
for $k=10$ with $n_{10}=29$. Therefore, we have $h(n,10)\leq\lceil\frac{3n}{5}\rceil-2$ for $n\geq 29$.
\bpr
By item $(iii)$, we have $|{\cal S}_{n,k}|\leq \lceil\frac{n}{2}\rceil+ \lceil\frac{n}{k}\rceil-\lceil\frac{n_k}{k}\rceil$ for $n_k\leq n<n_k+4k$.
A solution ${\cal S}_{n,k}$ satisfying $|{\cal S}_{n,k}|\leq \lceil\frac{n}{2}\rceil+ \lceil\frac{n}{k}\rceil-\lceil\frac{n_k}{k}\rceil$ 
for $n\geq n_k$ can be constructed inductively as follows.
Let $n_k+4ik\leq n<4(i+1)k+n_k$ and 
assume that $|{\cal S}_{n,k}|\leq \lceil\frac{n}{2}\rceil+ \lceil\frac{n}{k}\rceil-\lceil\frac{n_k}{k}\rceil$ for $n_k+4(i-1)k\leq n<4ik+n_k$.
Use the solution ${\cal S}_{4k,k}$ to sort to the left the first $4k$ pegs while ignoring the remaining  $n-4k$ pegs.
Then sort the ignored $n-4k$ pegs using the solution ${\cal S}_{n-4k,k}$. Complete the solution ${\cal S}_{n,k}$ by the 
elementary $4$ moves which append the sorted $4k$ pegs to the sorted $n-4k$ pegs. Items $(ii)$ and $(iv)$ guarantee the validity
of this solution ${\cal S}_{n,k}$, which takes at most $|{\cal S}_{n-4k,k}|+|{\cal S}_{4k,k}|+4$ moves; that is, by the induction hypothesis 
and item $(iii)$, at most
$\lceil\frac{n-4k}{2}\rceil+ \lceil\frac{n-4k}{k}\rceil-\lceil\frac{n_k}{k}\rceil+2k+4=\lceil\frac{n}{2}\rceil+ \lceil\frac{n}{k}\rceil-\lceil\frac{n_k}{k}\rceil$
moves.
\epr
\ep

\subsection{Proof of Proposition~\ref{3moves} }\label{proofs-3}
We exhibit solutions ${\cal S}_{n,3}$ in $\lceil\frac{n}{2}\rceil$ moves
for  $n\not\equiv 0\pmod{4}$  and in $\lceil\frac{n}{2}\rceil+1$ moves for $n\equiv 0\pmod{4}$.

\subsubsection{Case $n\equiv 1\pmod{4}$}\label{bbase}
\begin{figure}[htb]
\begin{center}
\epsfig{file=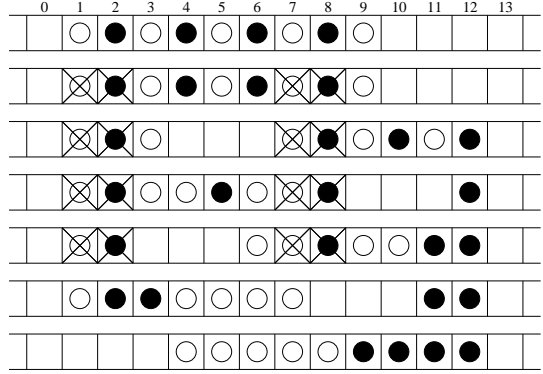,width=7.1cm}
\caption{Sorting $9$ pegs using the solution for $5$}
\label{Berge9-3}
\end{center}
\end{figure}

We have ${\cal S}_{5,3} = \{\ 6\ 2\ 5\ 1\ \}$ and ${\cal S}_{n,3}$ can be constructed inductively as follows.
Let $n=4i+1\geq 9$ and assume we have a solution ${\cal S}_{4i-3,3}$ taking $\lceil\frac{4i-3}{2}\rceil$ moves.
First ignore the $4$ pegs in positions $1,2,2i+3$
and $2i+4$ and sort the remaining $4i-3$ pegs using the solution ${\cal S}_{4i-3,3}$. Then complete the solution ${\cal S}_{4i+1,3}$ by 
the $2$ moves $\{ \ 3\ \ 2i + 4 \ \ 1 \ \}$.
The solution ${\cal S}_{4i+1,3}$ takes $\lceil\frac{4i-3}{2}\rceil+2=\lceil\frac{n}{2}\rceil$ moves.
Note that the solution ${\cal S}_{4i-3,3}$ can be performed while ignoring the
$4$ pegs in positions $1,2,2i+3$ and $2i+4$ because these pegs are not moved as, by induction, 
the solution $S_{4i+1,3}$ does not include among its entries any of $-1$, 0, $2i+1$, or $2i+2$ in the first $2i+1$ moves for $i\geq 1$.
More precisely, with ${\cal S}_{n,3}^j$ denoting the $j$-th entry of the solution ${\cal S}_{n,3}$, we have:
\begin{displaymath}
{\cal S}_{4i+1,3}^j= \left\{
\begin{array}{ll}
{\cal S}_{4i-3,3}^j + 2 & \mbox{ for } 1\leq {\cal S}_{4i-3,3}^j \leq 2i-2\\
{\cal S}_{4i-3,3}^j+ 4 & \mbox { for }  2i+1\leq  {\cal S}_{4i-3,3}^j   
\end{array}
\right.
\end{displaymath}
See Table~\ref{Bergen-3_1} for the first solutions ${\cal S}_{n,3}$ for $n=5,9,13$ and $17$
and Figure~\ref{Berge9-3} illustrating the induction from ${\cal S}_{5,3}$ to ${\cal S}_{9,3}$.

\begin{table}[htb]
\caption{First solutions for sorting $n$ pegs in $\lceil\frac{n}{2}\rceil$ Berge $3$-moves for $n\equiv 1\pmod 4$}
\label{Bergen-3_1}
\begin{center}
\begin{tabular}{r c l l l l l l l l l l l l l l l l l l l}
${\cal S}_{5,3}$  & =  & \{  6 \  2 \   5 \  1  \} \\
${\cal S}_{9,3}$ & =  & \{   10 \  4 \  9 \  3 \ 8 \  1  \} \\
${\cal S}_{13,3}$ & =  & \{  14 \  6 \  13 \ 5 \ 12 \  3 \ 10 \  1  \} \\
${\cal S}_{17,3}$ & =  & \{  18 \  8 \  17 \ 7 \  16 \   5 \ 14 \  3 \ 12 \  1  \}  
\end{tabular}
\end{center}
\end{table}

\newpage
\subsubsection{Case $n\equiv 2\pmod{4}$}\label{base}
We have ${\cal S}_{6,3} = \{\ 7\ 2\ 6\ 1\ \}$ and ${\cal S}_{n,3}$ can be constructed inductively as follows.
Let $n=4i+2\geq 10$ and assume we have a solution ${\cal S}_{4i-2,3}$ taking $\lceil\frac{4i-2}{2}\rceil$ moves.
First ignore the $4$ pegs in positions $1,2,2i+3$
and $2i+4$ and sort the remaining $4i-2$ pegs using the solution ${\cal S}_{4i-2,3}$. Then complete the solution ${\cal S}_{4i+2,3}$ by 
the $2$ moves $\{ \ 3\ \ 2i + 4 \ \ 1 \ \}$. 
The solution ${\cal S}_{4i+2,3}$  takes $\lceil\frac{4i-2}{2}\rceil+2=\lceil\frac{n}{2}\rceil$ moves.
Note that the solution ${\cal S}_{4i-2,3}$ can be performed while ignoring the
$4$ pegs in positions $1,2,2i+3$ and $2i+4$ because, by an argument similar to the one used in 
Section~\ref{bbase},  these pegs are not moved.
See Table~\ref{Bergen-3_2} for the first solutions ${\cal S}_{n,3}$ for  $n=6,10,14$ and $18$.

\begin{table}[htb]
\caption{First solutions for sorting $n$ pegs in $\lceil\frac{n}{2}\rceil$ Berge $3$-moves for $n\equiv 2\pmod 4$}
\label{Bergen-3_2}
\begin{center}
\begin{tabular}{r c l l l l l l l l l l l l l l l l l l l}
${\cal S}_{6,3}$ & = & \{  7 \ 2 \ 6 \ 1  \}\\
${\cal S}_{10,3}$ & = & \{  11 \ 4 \ 10 \ 3 \ 8 \ 1  \}\\
${\cal S}_{14,3}$ & = & \{  15 \ 6 \ 14 \ 5 \ 12 \ 3 \ 10 \ 1  \}\\
${\cal S}_{18,3}$ & = & \{  19 \ 8 \ 18 \ 7 \ 16 \ 5 \ 14 \ 3 \ 12 \ 1  \}
\end{tabular}
\end{center}
\end{table}

The following lemma can be easily checked by induction.
\bl{lemma_a}$\;$
\begin{itemize}
\item[$(i)$]
For $n\equiv 2\pmod 4$, the solutions ${\cal S}_{n,3}$ shift the string three spaces to the right overall.
\item[$(ii)$]
For $n\equiv 2\pmod 4$, the solutions ${\cal S}_{n,3}$ place the $\lceil\frac{n}{2}\rceil$ white pegs to the left of the $\lfloor\frac{n}{2}\rfloor$ 
black pegs\end{itemize}
\el

\subsubsection{Case $n\equiv 3\pmod{4}$}
We have ${\cal S}_{7,3}=\{ \ -\!2 \ \ 4 \ -\!1 \ \ 3 \ -\!2 \ \}$.  
Let $n=4i+3\geq 11$,  first perform the move $\{\ -2\ \ 4i\ \}$. Then, ignore the peg at position $4i+3$ and 
sort the remaining $4i+2$ pegs  using the solution ${\cal S}_{4i+2,3}$, see Section~\ref{base}.
Lemma~\ref{lemma_a} guarantees the validity of this solution ${\cal S}_{4i+3,3}$ which takes 
$\lceil \frac{4i+2}{2}\rceil+1=\lceil \frac{n}{2}\rceil$ moves. See Table~\ref{Bergen-3_3} for the first 
solutions ${\cal S}_{n,3}$ for $n=7,11,15$ and $19$.

\begin{table}[htb]
\caption{First solutions for sorting $n$ pegs in $\lceil\frac{n}{2}\rceil$ Berge $3$-moves for $n\equiv 3\pmod 4$}
\label{Bergen-3_3}
\begin{center}
\begin{tabular}{r c l l l l l l l l l l l l l l l l l l l l l l}
${\cal S}_{7,3}$& = & \{  -2 \ 4 \ -1 \ 3 \ -2  \}\\
${\cal S}_{11,3}$ & = & \{  -2 \ 8 \ 1 \ 7 \ 0 \ 5 \ -2  \}\\
${\cal S}_{15,3}$ & = & \{  -2 \ 12 \ 3 \ 11 \ 2 \ 9 \ 0 \ 7 \ -2  \}\\
${\cal S}_{19,3}$& = & \{  -2 \ 16 \ 5 \ 15 \ 4 \ 13 \ 2 \ 11 \ 0 \ 9 \ -2  \}
\end{tabular}
\end{center}
\end{table}\

\subsubsection{Case $n\equiv 0\pmod{4}$}
Although we found solutions in $\lceil\frac{n}{2}\rceil$ moves
for $n \equiv 0\pmod{4}$, $20 \leq n \leq 50$, we could not find solutions
in $\lceil\frac{n}{2}\rceil$ moves for all $n$.
However, solutions ${\cal {\bar S}}_{4i,3}$ in $\lceil\frac{n}{2}\rceil+1$ moves
can be constructed as follows. Let $n=4i\geq 16$, first perform the 2 moves $\{\   4i + 1 \ \ 2 \  \ 4i - 3 \ \}$.
Then, ignore the six leftmost pegs, and the four rightmost pegs and sort  the remaining $4i-10$ pegs using the solution ${\cal S}_{4i-10,3}$ shifted six spaces to the right, see Section~\ref{base}.
Finally, perform the 4 moves  $\{\ 7 \ \ 4i \ \ 6 \ \  2i + 2 \  \ 1\ \}$ to complete the solution ${\cal {\bar S}}_{4i,3}$.
Lemma~\ref{lemma_a} guarantees the validity of this solution ${\cal {\bar S}}_{4i,3}$ which takes 
$2+\lceil \frac{4i-10}{2}\rceil + 4 = \lceil \frac{n}{2}\rceil + 1$ moves.
See Table~\ref{Bergen-3_4} for the first solutions ${\cal {\bar S}}_{n,3}$ for
$n=16,20$ and $24$.

\begin{table}[htb]
\caption{First solutions for sorting $n$ pegs in $\lceil\frac{n}{2}\rceil+1$ Berge $3$-moves for $n\equiv 0\pmod 4$}
\label{Bergen-3_4}
\begin{center}
\begin{tabular}{r c l l l l l l l l l l l l l l l l l l l l l l}
${\cal {\bar S}}_{16,3}$ & = & \{ 17 \ 2 \ 13 \ 8 \ 12 \ 7 \ 16 \ 6 \ 10 \ 1 \}\\
${\cal {\bar S}}_{20,3}$ & = & \{ 21 \ 2 \ 17 \ 10 \ 16 \ 9 \ 14 \ 7 \ 20 \ 6 \ 12 \ 1 \}\\
${\cal {\bar S}}_{24,3}$ & = & \{ 25 \ 2 \ 21 \ 12 \ 20 \ 11 \ 18 \ 9 \ 16 \ 7 \ 24 \ 6 \ 14 \ 1 \}
\end{tabular}
\end{center}
\end{table}
\noindent
While we could not exhibit solutions in $\lceil\frac{n}{2}\rceil$ moves for all $n\equiv 0\pmod 4$,
we believe that such solutions exist for $n\geq 20$, i.e., the proposed solutions ${\cal {\bar S}}_{4i,3}$ are not optimal,
except for ${\cal {\bar S}}_{16,3}$. See Table~\ref{Bergen-3_4bis} for 
optimal solutions in $\lceil\frac{n}{2}\rceil+1$ moves for $n=12$ and $16$, and Table~\ref{Bergen-3_4ter} for optimal
solutions in $\lceil\frac{n}{2}\rceil$ moves  for $n=8,20,24,28$ and $32$.

\begin{table}[htb]
\caption{Solutions for sorting $n$ pegs in $\lceil\frac{n}{2}\rceil+1$ Berge $3$-moves for $n=12$ and $16$}
\label{Bergen-3_4bis}
\begin{center}
\begin{tabular}{r c l l l l l l l l l l l l l l l l l l l l l l}
${\cal S}_{12,3}$ & = & \{  13 \ 2 \ 5 \ 11 \ 3 \ 12 \ 6 \ 1  \}\\
${\cal {\bar S}}_{16,3}$ & = & \{ 17 \ 2 \ 13 \ 8 \ 12 \ 7 \ 16 \ 6 \ 10 \ 1 \}\\
\end{tabular}
\end{center}
\end{table}

\begin{table}[htb]
\caption{Solutions for sorting $n$ pegs in $\lceil\frac{n}{2}\rceil$ Berge $3$-moves for $n\!=\!8,20,24,28$ and $32$}
\label{Bergen-3_4ter}
\begin{center}
\begin{tabular}{r c l l l l l l l l l l l l l l l l l l l l l l}
${\cal S}_{8,3}$& = & \{ 9 \ 2 \ 7 \ 3 \ 9 \}\\
${\cal S}_{20,3}$& = & \{ 21 \ 2 \ 7 \ 12 \ 17 \} $\cup$ \{ 24 \ 13 \ 22 \ 6 \ 1 \} $\cup$ \{ 17 \ 8 \ 24 \}\\
${\cal S}_{24,3}$& = & \{ 25 \ 6 \ 13 \ 18 \} $\cup$ \{ -2 \ 4 \ 8 \ 24 \ 14 \ 22 \} $\cup$ \{ 18 \ 3 \ 12 \ -1 \ 25 \}\\
${\cal S}_{28,3}$& = & \{ 29 \ 2 \ 7 \ 16 \ 23 \ 12 \} $\cup$ \{ 32 \ 17 \ 30 \ 25 \ 21 \ 6 \ 1 \} $\cup$ \{ 12 \ 23 \ 8 \ 32 \}\\
${\cal S}_{32,3}$& = & \{ 33 \ 2 \ 7 \ 12 \ 17 \ 24 \} $\cup$ \{ 36 \ 6 \ 31 \ 13 \ 29 \ 19 \ 1 \} $\cup$ \{ 24 \ 11 \ 35 \ 18 \ 28 \ 4 \}
\end{tabular}
\end{center}
\end{table}

 \section{Related Questions}\label{generalisation}
Other extensions of Berge's original questions include sorting any $n$ string:
\begin{enumerate}
\item[$(a_1)$] 
Besides the alternating string,  which other string requires exactly $h(n,k)$  Berge $k$-moves? 
\item[$(a_2)$] 
What is the minimum number of Berge $k$-moves required to sort any $n$ string? 
\item [$(a_3)$] 
Given a pair of strings, can we rearrange one into the other by Berge $k$-moves?
\end{enumerate}
Associating the white and black colors to  $0$ and $1$, the original $\{0,1\}$-valued string could be generalized to 
$\{0,1,\dots,m\}$-valued strings where $m$ is the number of colors; the final string being 
$0\dots 0\: 1\dots 1\dots m\dots m$:
\begin{enumerate}
\item[$(b_1)$] 
What is the minimum number of Berge $k$-moves required to sort a string consisting of $m$ different integers
- each integer being represented by the same number of pegs? 
\item [$(b_2)$] 
In particular, what is the minimum number
of Berge $k$-moves required to sort a string consisting of $n$ different integers.
\end{enumerate}
Generalizing to moves of $k$-by-$k$ blocks in the plane could also be considered.
 Similar questions were raised for $2$-moves in~\cite{AD0X}.\\\\

\noindent {\bf Acknowledgments}.
Research supported  by the 
Natural Sciences and Engineering Research Council of Canada
under the Canada Research Chair and the Discovery Grant  programs.

\begin{table}[htb]
\caption{Values of $h(n,k)-\lceil\frac{n}{2}\rceil$  for $k \leq 14$ and $n \leq 50$}\label{substance}
\begin{center}
\begin{small}
\begin{tabular}{|c|ccccccccccccc|c|}
\hline
$n\backslash k$ & 2 & 3 & 4 & 5 & 6 & 7 & 8 & 9 & 10 & 11 & 12 & 13 & 14 \\
\hline
5 & 0 & 0 & -- & -- & -- & -- & -- & -- & -- & -- & -- & -- & -- \\
6 & 0 & 0 & 3 & -- & -- & -- & -- & -- & -- & -- & -- & -- & -- \\
7 & 0 & 0 & 0 & 2 & -- & -- & -- & -- & -- & -- & -- & -- & -- \\
8 & 0 & 0 & 1 & 2 & 3 & -- & -- & -- & -- & -- & -- & -- & -- \\
9 & 0 & 0 & 0 & 1 & 2 & 3 & -- & -- & -- & -- & -- & -- & -- \\
10 & 0 & 0 & 1 & 1 & 1 & 3 & 6 & -- & -- & -- & -- & -- & -- \\
11 & 0 & 0 & 0 & 1 & 1 & 2 & 4 & 6 & -- & -- & -- & -- & -- \\
12 & 0 & 1 & 1 & 1 & 1 & 2 & 3 & 5 & 10 & -- & -- & -- & -- \\
13 & 0 & 0 & 0 & 1 & 1 & 1 & 2 & 3 & 4 & 11 & -- & -- & -- \\
14 & 0 & 0 & 0 & 1 & 2 & 2 & 2 & 2 & 4 & 6 & 15 & -- & -- \\
15 & 0 & 0 & 0 & 0 & 1 & 1 & 1 & 2 & 2 & 4 & 7 & 14 & -- \\
16 & 0 & 1 & 0 & 1 & 1 & 0 & 2 & 2 & 3 & 3 & 5 & 9 & 21 \\
17 & 0 & 0 & 0 & 0 & 0 & 1 & 1 & 1 & 2 & 2 & 3 & 5 & 9 \\
18 & 0 & 0 & 0 & 1 & 1 & 1 & 1 & 1 & 2 & 3 & 3 & 4 & 7 \\
19 & 0 & 0 & 0 & 0 & 0 & 0 & 1 & 1 & 1 & 1 & 2 & 3 & 4 \\
20 & 0 & 0 & 0 & 0 & 1 & 1 & 1 & 1 & 2 & 2 & 3 & 3 & 4 \\
21 & 0 & 0 & 0 & 0 & 0 & 0 & 0 & 0 & 1 & 1 & 2 & 2 & 3 \\
22 & 0 & 0 & 0 & 0 & 1 & 1 & 1 & 1 & 1 & 2 & 1 & 2 & 3 \\
23 & 0 & 0 & 0 & 0 & 0 & 0 & 0 & 0 & 1 & 1 & 1 & 1 & 2 \\
24 & 0 & 0 & 0 & 0 & 0 & 0 & 1 & 1 & 1 & 1 & 1 & 1 & 2 \\
25 & 0 & 0 & 0 & 0 & 0 & 0 & 0 & 0 & 0 & 0 & 1 & 1 & 2 \\
26 & 0 & 0 & 0 & 0 & 0 & 0 & 1 & 0 & 0 & 1 & 1 & 1 & 2 \\
27 & 0 & 0 & 0 & 0 & 0 & 0 & 0 & 0 & 0 & 0 & 1 & 1 & 1 \\
28 & 0 & 0 & 0 & 0 & 0 & 0 & 0 & 0 & 1 & 1 & 1 & 1 & 1 \\
29 & 0 & 0 & 0 & 0 & 0 & 0 & 0 & 0 & 0 & 0 & 0 & 1 & 0 \\
30 & 0 & 0 & 0 & 0 & 0 & 0 & 0 & 0 & 0 & 1 & 1 & 1 & 1 \\
31 & 0 & 0 & 0 & 0 & 0 & 0 & 0 & 0 & 0 & 0 & 0 & 0 & 0 \\
32 & 0 & 0 & 0 & 0 & 0 & 0 & 0 & 0 & 0 & 0 & 0 & 0 & 1 \\
33 & 0 & 0 & 0 & 0 & 0 & 0 & 0 & 0 & 0 & 0 & 0 & 0 & 0 \\
34 & 0 & 0 & 0 & 0 & 0 & 0 & 0 & 0 & 0 & 0 & 0 & 0 & 1 \\
35 & 0 & 0 & 0 & 0 & 0 & 0 & 0 & 0 & 0 & 0 & 0 & 0 & 0 \\
36 & 0 & 0 & 0 & 0 & 0 & 0 & 0 & 0 & 0 & 0 & 0 & 0 & 1 \\
37 & 0 & 0 & 0 & 0 & 0 & 0 & 0 & 0 & 0 & 0 & 0 & 0 & 0 \\
38 & 0 & 0 & 0 & 0 & 0 & 0 & 0 & 0 & 0 & 0 & 0 & 0 & 0 \\
39 & 0 & 0 & 0 & 0 & 0 & 0 & 0 & 0 & 0 & 0 & 0 & 0 & 0 \\
40 & 0 & 0 & 0 & 0 & 0 & 0 & 0 & 0 & 0 & 0 & 0 & 0 & 0 \\
41 & 0 & 0 & 0 & 0 & 0 & 0 & 0 & 0 & 0 & 0 & 0 & 0 & 0 \\
42 & 0 & 0 & 0 & 0 & 0 & 0 & 0 & 0 & 0 & 0 & 0 & 0 & 0 \\
43 & 0 & 0 & 0 & 0 & 0 & 0 & 0 & 0 & 0 & 0 & 0 & 0 & 0 \\
44 & 0 & 0 & 0 & 0 & 0 & 0 & 0 & 0 & 0 & 0 & 0 & 0 & 0 \\
45 & 0 & 0 & 0 & 0 & 0 & 0 & 0 & 0 & 0 & 0 & 0 & 0 & 0 \\
46 & 0 & 0 & 0 & 0 & 0 & 0 & 0 & 0 & 0 & 0 & 0 & 0 & 0 \\
47 & 0 & 0 & 0 & 0 & 0 & 0 & 0 & 0 & 0 & 0 & 0 & 0 & 0 \\
48 & 0 & 0 & 0 & 0 & 0 & 0 & 0 & 0 & 0 & 0 & 0 & 0 & 0 \\
49 & 0 & 0 & 0 & 0 & 0 & 0 & 0 & 0 & 0 & 0 & 0 & 0 & 0 \\
50 & 0 & 0 & 0 & 0 & 0 & 0 & 0 & 0 & 0 & 0 & 0 & 0 & 0 \\
\hline
\end{tabular}
\end{small}
\end{center}
\end{table}

\noindent
{\small Antoine Deza, William Hua}\\
Advanced Optimization Laboratory,
Department of Computing and Software,\\
McMaster University, Hamilton, Ontario, Canada. \\
{\em Email}: deza, huaw{\small @}mcmaster.ca.

\end{document}